\documentclass[10pt,leqno]{article}

\usepackage{amsmath,amssymb,amsthm,mathrsfs,dsfont}

\usepackage[margin=3cm]{geometry} 

\usepackage{titlesec,hyperref}

\usepackage{pgfplots}
\usepackage{color}

\usepackage{fancyhdr}
\titleformat{\subsection}{\it}{\thesubsection.\enspace}{1pt}{}

\newtheorem{theo}{Theorem}[section]
\newtheorem{lemm}[theo]{Lemma}
\newtheorem{defi}[theo]{Definition}

\newtheorem{rema}[theo]{Remark}
\numberwithin{equation}{section}

\allowdisplaybreaks 

\newcommand\ep{{\varepsilon}} 
\newcommand{\num}{pi}

\begin{document}
\title{Local bifurcation of steady almost periodic water waves with constant vorticity
\hspace{-4mm}
}

\author{Wei $\mbox{Luo}^1$\footnote{E-mail:  luowei23@mail2.sysu.edu.cn} \quad and\quad
 Zhaoyang $\mbox{Yin}^{1,2}$\footnote{E-mail: mcsyzy@mail.sysu.edu.cn}\\
 $^1\mbox{Department}$ of Mathematics,
Sun Yat-sen University,\\ Guangzhou, 510275, China\\
$^2\mbox{Faculty}$ of Information Technology,\\ Macau University of Science and Technology, Macau, China}
\date{}
\maketitle
\hrule

\begin{abstract}
In this paper we mainly investigate the traveling wave solution of the two dimensional Euler equations with gravity at the free surface over a flat bed. We assume that the free surface is almost periodic in the horizontal direction. Using conformal mappings, one can change the free boundary problem into a fixed boundary problem with some unknown functions in the boundary condition. By virtue of the Hilbert transform, the problem is equivalent to a quasilinear pseudodifferential equation for a almost periodic function of one variable. The bifurcation theory ensures us to obtain a existence result. Our existence result generalizes and covers the recent result in \cite{Constantin2011v}. Moreover, our result implies a non-uniqueness result at the same bifurcation point.\\

\vspace*{5pt}

\vspace*{5pt}
\end{abstract}

\vspace*{10pt}

\tableofcontents
\section{Introduction}
  In this paper we consider the following free-boundary problem \cite{Constantin2004}:
  \begin{align}\label{1}
\left\{
\begin{array}{ll}
\Delta \psi=-\gamma,  \quad \text{in}\quad \Omega, \\[1ex]
\psi=-m,     \quad  x\in\mathbb{R}, ~~y=0, \\[1ex]
\psi=0, \quad    x\in\mathbb{R}, ~~ y=\eta(x),\\[1ex]
|\nabla \psi|^2+2gY=Q   \quad x\in\mathbb{R},  y=\eta(x).\\[1ex]
\end{array}
\right.
\end{align}
The above problem can be deduced from the travelling wave solution of the Euler equations with constant vorticity (See \cite{Constantin2004,Wahlen}). Here $\psi$ represents the stream function. The velocity field $v=(\psi_Y,-\psi_X)$ in a frame moving at the constant wave speed. The vorticity $w=\partial_x v^2-\partial_y v^1=\gamma$, where $\gamma$ is a given constant. The constant $g$ is the gravitational acceleration, while the constant $m$ is the relative mass flux. The constant $Q$ is deduced from Bernoulli's law which is related to the hydraulic head (See \cite{Constantin2004}). The domain $\Omega$ belongs to a two dimensional $(X,Y)-$plane which is given by
\begin{align}
\Omega=\{(x,y)\in\mathbb{R}^2|0<y<\eta(x)\},
\end{align}
where $\eta(x)$ is an unknown curve representing the free surface of the water. In this paper, we will assume that $\eta(x)$ is an almost periodic function.

The irrotational $(\gamma=0)$ periodic travelling water waves have been studied in a long time. The first original study was performed by G. Stokes \cite{Stokes}. For irrotational flows, a classical approach is to use the hodograph transformation mapping the unknown domain into a fixed strip. More precisely, taking advantage of the changing variables:
\[(X,Y)\mapsto(X,\frac{Y}{\eta(x)}),\]
one can reformulate $(\ref{1})$ as a one dimensional problem formed by a nonlinear singular integral equation. And then using the bifurcation theory, one can obtain the existence theory and geometric properties of the water surface (See \cite{Keady,Toland,Groves}).

A more recent approach is to use the conformal mapping:
 \[(X,Y)\mapsto (U(X,Y), V(X,Y)),\]
to change the domain $\Omega$ into a fixed strip. The problem (\ref{1}) is equivalent to solve a harmonic function in a strip with two different unknown boundary conditions. In order to solve the harmonic function, these two boundary conditions must be compatible. Taking advantage of the periodic Dirichlet-Neumann operator, one can change the problem $(\ref{1})$ into a one dimensional pseudodifferential equation (See \cite{Buffoni,Buffoni2,Toland2}). This approach can be extended to investigate the rotational water waves.

The periodic travelling water waves with nonzero vorticity have been studied widely in recent years. The existence, uniqueness, regularity and geometric properties of solutions to $(\ref{1})$ were recently investigated by many researchers. Several existence results were proved under different conditions by using different methods (\cite{Constantin2004,Constantin2011,Constantin2016,Wahlen,Groves2008}). The symmetry properties of the water surface was studied in \cite{Constantin2007.Duke,Constantin2004.EJAM,Constantin2004.JFM}. The uniqueness results were obtained in \cite{Ehrnstrom,Ehrnstrom2005}. The geometric properties and regularity of solutions have been studied in \cite{Constantin2006,Constantin2007.CPAM,Constantin2011.Escher}.

In the previous works, the periodic water waves are focused on. However, in a real world, some water waves may not be periodic, but almost periodic. To our best knowledge, the almost periodic water wave problem has not been studied yet. In this paper, we will study this problem.  Being inspired by the recent work of A. Constantin et al in \cite{Constantin2011v,Constantin2016}, we will use the conformal mapping to change $(\ref{1})$ into a pseudodifferential equation. Then we apply the local bifurcation theory to construct a family of solutions that are a small perturbations of laminar flows. The main difficulty is to choose two suitable Banach spaces $X,Y$ composed of almost periodic functions such that the bifurcation theory can be applied in. Indeed, our obtained result generalizes and covers the recent result in \cite{Constantin2011v}. Moreover, we construct a solution that is not an even function, which is quite different from the previous existence results (See, \cite{Constantin2004,Constantin2011,Constantin2011v,Constantin2016,Wahlen,Groves2008}). Our result also implies that for any given constants $\gamma$, $m$, $Q$ and the conformal mean depth, the solution of $(\ref{1})$ is not unique at the same bifurcation point.

The remainder of the paper is organized as follows. In Section 2 we introduce some preliminaries about almost periodic functions and the Dirichlet-Neumman operator which will be used in sequel. In Section 3 we reformulate the free boundary problem (\ref{1}) to a pseudodifferential equation. Section 4 is devoted to studying the existence result of the new reformulated problem by using the local bifurcation theory. Moreover, we give an example that is quiet different from that of \cite{Constantin2011v,Constantin2016}.

\section{Besicovitch's almost periodic functions}
In this section, we recall some basic properties for almost periodic functions introduced by A. S. Besicovitch.
\begin{defi}\cite{Besicovitch}
Let $u:\mathbb{R}\rightarrow \mathbb{R}$ be a smooth function. $u$ is called an almost periodic function if and only if there exist two real sequences $(\alpha_k)_{k\geq 0}$ and $(\beta_k)_{k\geq0}$ such that
\[u=\overline{u}+\sum_{k=0}a_k\cos(\alpha_k x)+\sum_{k=0}b_k\sin(\beta_k x),\]
where $\overline{u}=\lim_{x\rightarrow\infty}\frac{1}{2x}\int^x_{-x}u(y)dy$, and  $$a_k=\lim_{x\rightarrow\infty}\frac{1}{x}\int^x_{-x}u(y)\cos(\alpha_k y) dy, ~~~b_k=\lim_{x\rightarrow\infty}\frac{1}{x}\int^x_{-x}u(y)\sin(\beta_k y) dy.$$
\end{defi}
\begin{rema}
$(i)$ Without loss of generality, we may assume that $\alpha_k>\alpha_j$ and $\beta_k>\beta_j$ if $k>j$. \\
$(ii)$ A periodic function is also an almost periodic function with $\alpha_k=\beta_k=k$.
\end{rema}

\begin{defi}\cite{Besicovitch}
Let $p\in[1,+\infty)$. The Besicovitch space $\mathcal{B}^p(\mathbb{R})$ consists of all the functions $u\in L^p_{loc}(\mathbb{R})$ such that
\[\|u\|_{\mathcal{B}^p(\mathbb{R})}\triangleq \lim_{x\rightarrow\infty} \bigg(\frac{1}{x}\displaystyle\int^x_{-x}|u|^pdy\bigg)^\frac{1}{p}.\]
\end{defi}

\begin{rema}\cite{Besicovitch}
Suppose that $f$ is a bounded function with compact support, then $\|f\|_{\mathcal{B}^p(\mathbb{R})}=0$. Thus, $\|\cdot\|_{\mathcal{B}^p(\mathbb{R})}$ is a semi-norm. In order to get a Banach space, one has to quotient out by these functions.
\end{rema}

\begin{lemm}\cite{Besicovitch}
Let $L_0$ be the class consists of all functions $\|u\|_{\mathcal{B}^p(\mathbb{R})}=0$. Then, the quotient Besicovitch space $\mathcal{B}^p(\mathbb{R})/L_0$ is a Banach space. Moreover, if $u\in \mathcal{B}^2(\mathbb{R})/L_0$ then there exist two real sequences $(\alpha_k)_{k\geq0}$ and $(\beta_k)_{k\geq0}$ such that
\[u=\overline{u}+\sum_{k=0}a_k\cos(\alpha_k x)+\sum_{k=0}b_k\sin(\beta_k x)\quad \text{in} \quad \mathcal{B}^2(\mathbb{R}),\]
\[\|u\|_{\mathcal{B}^2(\mathbb{R})}\approx \overline{u}^2+(\sum_{k>0}(a^2_k+b^2_k))^{\frac{1}{2}}.\]
\end{lemm}

\begin{rema}\cite{Besicovitch}
From the above lemma, one can extend the definition of the almost periodic function on the Besicovitch space.
\end{rema}

\begin{rema}\cite{Besicovitch}
$\mathcal{B}^2(\mathbb{R})/L_0$ is a Hilbert space equipped with inner product
$$<u,v>\triangleq \lim_{x\rightarrow\infty}\frac{1}{x}\int^x_{-x}u(y)v(y)dy.$$
Moreover, $\{1, \cos(\alpha_k x), \sin(\beta_k x)|k\geq 0\}$ is a standard orthogonal basis.
\end{rema}

For any $h>0$, let $\mathcal{R}_h$ be the strip
\[\mathcal{R}_h=\{(x,y)\in\mathbb{R}^2:-h<y<0\}.\]
For any $n\in \mathbb{N}^{+}$ and $\alpha\in(0,1)$, we denote by $C^{n,\alpha}$ the standard H\"{o}lder space. For any $w\in C^{n,\alpha}(\mathbb{R})\cap \mathcal{B}^2(\mathbb{R})/L_0$, let $W\in C^{n,\alpha}(\overline{\mathcal{R}}_h)$ be the unique solution of
\begin{equation}
\left\{
\begin{array}{ll}
\Delta W=0, \quad \text{in} \quad \mathcal{R}_h, \\[1ex]
W(x,-h)=0,  \quad  x\in\mathbb{R}, \\[1ex]
W(x,0)=w(x), \quad x\in\mathbb{R}.
\end{array}
\right.
\end{equation}
For $n\geq 1$, we define that
\[\mathcal{G}_h(w)(x)=W_y(x,0), \quad x\in\mathbb{R}. \]
The mapping $w\mapsto\mathcal{G}_h(w)$ is called the Drichlet-Neumann operator for a strip. If $c$ is a constant and $w\equiv c$, then
\[W(x,y)=\frac{c}{h}(y+h), \quad (x,y)\in \mathcal{R}_h, \]
and hence
\[\mathcal{G}_h(c)=\frac{c}{h}.\]
Since $w\in \mathcal{B}^2(\mathbb{R})/L_0$, it follows that
\[w=\overline{w}+\sum_{k=0}a_k\cos(\alpha_k x)+\sum_{k=0}b_k\sin(\beta_k x), \quad x\in \mathbb{R}, \]
and then
\[W(x,y)=\frac{\overline{w}}{h}(y+h)+\sum^\infty_{k=0}a_k\frac{\sinh(\alpha_k(y+h))}{\sinh(\alpha_k h)}\cos(\alpha_k x)+\sum^\infty_{k=0}b_k\frac{\sinh(\beta_k(y+h))}{\sinh(\beta_k h)}\sin(\beta_k x),\]
which implies that
\begin{align}\label{2.1}
\mathcal{G}_h(w)=\frac{\overline{w}}{h}+\sum^\infty_{k=0}\alpha_k a_k\coth(\alpha_k h)\cos(\alpha_k x)+\sum^\infty_{k=0}\beta_k a_k\coth(\beta_k h)\cos(\beta_k x).
\end{align}
From the above equality, we see that the mapping $w\mapsto\mathcal{G}_h(w)$ is a bounded linear operator from $C^{n,\alpha}(\mathbb{R})\cap \mathcal{B}^2(\mathbb{R})/L_0$ to $C^{n-1,\alpha}(\mathbb{R})\cap \mathcal{B}^2(\mathbb{R})/L_0$.

\section{Reformulation of the free-boundary problem}
In this section we present the reformulation of the free-boundary problem $(1.1)$. We are interested in the solution $(\eta, \psi)$ of the water-wave problem $(1.1)$ where $\eta \in C^{1,\alpha}(\mathbb{R})\cap \mathcal{B}^2(\mathbb{R})/L_0$ and $\psi\in C^{1,\alpha}(\overline{\Omega}) $ for some $\alpha\in(0,1)$. The main result of this section is that the free-boundary problem $(1.1)$ is equivalent to a quasi-linear pseudodifferential equation:
\begin{equation}\label{3.1}
\left\{
\begin{array}{ll}
\bigg\{\frac{m}{h}+\gamma\bigg(\mathcal{G}_h(\frac{\eta^2}{2})-\eta\mathcal{G}_h(\eta)\bigg)\bigg\}^2=(Q-2g\eta)(\eta'^2+\mathcal{G}_h(\eta)^2), \\[1ex]
\overline{\eta}=h, \quad \eta>0, \quad \text{for all}  \quad x\in\mathbb{R}.
\end{array}
\right.
\end{equation}

\begin{theo}
 Let $(\eta,\psi)$ be a solution of $(1.1)$ where $\eta \in C^{1,\alpha}(\mathbb{R})\cap \mathcal{B}^2(\mathbb{R})/L_0$ and $\psi\in C^{1,\alpha}(\overline{\Omega})$ with
 \[\Omega=\{(x,y)\in\mathbb{R}^2|0<y<\eta(x)\}.\]
 Then $\eta$ is a solution of the equation (\ref{3.1}). Conversely, if $\eta\in C^{1,\alpha}(\mathbb{R})\cap \mathcal{B}^2(\mathbb{R})/L_0 $ is a solution of (\ref{3.1}). Then there exists a function $\psi\in C^{1,\alpha}(\overline{\Omega})$ which is a solution of $(1.1)$.
 \begin{proof}
 Let $(\eta,\psi)$ be a solution of $(1.1)$. As $\eta\in C^{1,\alpha}(\mathbb{R})\cap \mathcal{B}^2(\mathbb{R})/L_0$, we can define that \[h\triangleq\overline{\eta}=\lim_{x\rightarrow\infty}\frac{1}{x}\int^x_{-x}\eta(y)dy.\]
 Consider the conformal mapping $U+iV:\mathcal{R}_h\rightarrow \Omega$ satisfies the following equation
  \begin{equation}\label{3.2}
\left\{
\begin{array}{ll}
\Delta V=0, \quad \text{in} \quad \mathcal{R}_h, \\[1ex]
V(x,-h)=0,  \quad  x\in\mathbb{R}, \\[1ex]
V(x,0)=\eta (x), \quad x\in\mathbb{R}.
\end{array}
\right.
\end{equation}

Denote that $\xi(x,y)\triangleq\psi(U(x,y),V(x,y))$. By directly calculating, we see that
\[\xi_x=\psi_{X}U_x+\psi_YV_x,\quad \xi_y=\psi_XU_y+\psi_YV_y,\]
\[\xi_{xx}=\psi_{XX}U^2_x+2\psi_{XY}U_xV_x+\psi_{X}U_{xx}+\psi_{YY}V^2_x+\psi_YV_{xx},\]
\[\xi_{yy}=\psi_{XX}U^2_y+2\psi_{XY}U_yV_y+\psi_{X}U_{yy}+\psi_{YY}V^2_y+\psi_YV_{yy}.\]
Using the fact that $U_x=V_y$, $V_x=-U_y$ and $\Delta U=\Delta V=0$, we obtain
\[\Delta \xi=\Delta\psi(V^2_x+V^2_y).\]
Since $\Delta\psi=-\gamma$, it follows that
\[\Delta \xi=-\gamma(V^2_x+V^2_y).\]
Note that $\psi(x,0)=-m$ and $\psi(x,\eta(x))=0$. By virtue of the boundary condition of (\ref{3.2}), we have
\[\xi(x,-h)=\psi(U(x,-h),0)=-m,\quad \xi(x,0)=\psi(U(x,0),\eta(x))=0.\]
On the other hand, we see that
\[\xi^2_x+\xi^2_y=(\psi^2_X(U,V)+\psi^2_Y(U,V))(V^2_x+V^2_y).\]
Taking advantage of the last equation of $(1.1)$, we deduce that
\[\xi^2_x+\xi^2_y=(Q-2gV)(V^2_x+V^2_y) \quad \forall x\in\mathbb{R} \quad \text{and} \quad y=0. \]

Let $\zeta\triangleq \xi+m+\frac{\gamma}{2}V^2$. The above relations implies that
  \begin{equation}\label{3.3}
\left\{
\begin{array}{ll}
\Delta \zeta=0, \quad \text{in} \quad \mathcal{R}_h, \\[1ex]
\zeta(x,-h)=0,  \quad  x\in\mathbb{R}, \\[1ex]
\zeta(x,0)=m+\frac{\gamma}{2}\eta^2 (x),  \quad x\in\mathbb{R}, \\[1ex]
(\zeta_y-\gamma VV_y)^2=(Q-2gV)(V^2_x+V^2_y),  \quad \text{at}~~(x,0)\quad  \forall x\in\mathbb{R}.
\end{array}
\right.
\end{equation}

By virtue of the definition of the Dirichlet-Neumman operator, we have
\[\zeta_y(x,0)=\mathcal{G}_h(m+\frac{\gamma}{2}\eta^2)=\frac{m}{h}+\gamma\mathcal{G}_h(\frac{\eta^2}{2}),~~~~V_y(x,0)=\mathcal{G}_h(\eta).\]
Plugging the above equality into (\ref{3.3}) yields that
\[\bigg\{\frac{m}{h}+\gamma\bigg(\mathcal{G}_h(\frac{\eta^2}{2})-\eta\mathcal{G}_h(\eta)\bigg)\bigg\}^2=(Q-2g\eta)(\eta'^2+\mathcal{G}_h(\eta)^2).\]
So $\eta$ is a solution of the equation (\ref{3.1}).

Conversely, supposing that $\eta$ is a solution of the equation (\ref{3.1}), we can solve (\ref{3.2}) and (\ref{3.3}) to construct $U$, $V$ and $\zeta$. Then define that $\xi\triangleq\zeta-m-\frac{\gamma}{2}V^2$. By solving the linear system
\[\xi_x=\psi_{X}U_x+\psi_YV_x,\quad \xi_y=\psi_XU_y+\psi_YV_y,\]
one can obtain $\psi_X$ and $\psi_Y$ which lead to a function $\psi$. Finally, it is easy to check that $\psi$ is the solution of $(1.1)$.
 \end{proof}
\end{theo}

\section{Existence theory}
In this section we prove the existence of solutions of (\ref{3.1}) by using a local bifurcation theory. Letting $\eta=w+h$, we can rewrite (\ref{3.1}) into
\begin{equation}\label{4.1}
\left\{
\begin{array}{ll}
\bigg\{\frac{m}{h}+\gamma\bigg(\mathcal{G}_h(\frac{w^2}{2})-w-\frac{h}{2}-w\mathcal{G}_h(w)\bigg)\bigg\}^2=(Q-2gw-2gh)(w'^2+\mathcal{G}_h(w)^2+2\mathcal{G}_h(w)+1), \\[1ex]
\overline{w}=0, \quad w(x)>-h \quad \text{for all} \quad x\in\mathbb{R}.
\end{array}
\right.
\end{equation}

Note that $w=0$ is a solution of (\ref{4.1}) if and only if
\[Q-2gh=(\frac{m}{h}-\frac{\gamma h}{2})^2.\]
Set
\[\lambda=\frac{m}{h}-\frac{\gamma h}{2}, \mu=Q-2gh-\lambda^2.\]
Then the equation (\ref{4.1}) can be rewritten as
\begin{equation}\label{4.2}
\left\{
\begin{array}{ll}
\bigg\{\lambda+\gamma\bigg(\mathcal{G}_h(\frac{w^2}{2})-w-w\mathcal{G}_h(w)\bigg)\bigg\}^2=(\lambda^2+\mu-2gw)(w'^2+\mathcal{G}_h(w)^2+2\mathcal{G}_h(w)+1), \\[1ex]
\overline{w}=0, \quad w(x)>-h \quad \text{for all}  ~~~x\in\mathbb{R}.
\end{array}
\right.
\end{equation}

In order to prove the existence result, we will apply the following local bifurcation theorem which was proved by M. G. Crandall and P. H. Rabinowitz in \cite{Rabinowitz}.

{\bf Notation: Let $\mathcal{L}$ be a linear operator between two Banach spaces. We denote by $\mathcal{N}(\mathcal{L})$ and $\mathcal{R}(\mathcal{L})$ its null space and range, respectively.}

\begin{theo}(Local bifurcation theorem \cite{Rabinowitz})\label{th1}
 Let $X$ and $Y$ be Banach spaces, $I$ be an open interval in $\mathbb{R}$ containing $\lambda^*$, and $F\in C(I\otimes X, Y)$. Suppose that \\
 $(i)$ ~~$F(\lambda, 0)=0$ for all $\lambda\in I$; \\
 $(ii)$ ~~$\partial_\lambda F,\partial_u F$ and $\partial^2_{\lambda, u}F$ exist and are continuous;\\
 $(iii)$~~ $\mathcal{N}(\partial_uF(\lambda^*,0))$ and $Y/\mathcal{R}(\partial_uF(\lambda^*,0))$ are one-dimensional, with the null space generated by $u^*$;\\
 $(iv)$ the transversality condition condition $\partial^2_{\lambda^*, 0}F(\lambda^*,0)(1,u^*)\notin \mathcal{R}(\partial_uF(\lambda^*,0))$ holds.\\
  Then there exists a continuous local bifurcation curve $\{(\lambda(s),u(s)):|s|<\ep\}$ with $\ep>0$ sufficiently small such that $(\lambda(0),u(0))=(\lambda^*,0)$ and
  $$\{(\lambda.u)\in \mathcal{O}:u\neq0, F(\lambda,u)=0\}=\{(\lambda(s),u(s)):0<|s|<\ep\}$$
  for some neighborhood $\mathcal{O}$ of $(\lambda^*,0)\in I\times X$. Moreover, we have
  \[ u(s)=su^*+o(s) \quad \text{in}\quad X, \quad |s|<\ep,\]
  and if $\partial^2_{u}F$ is also continuous, the the curve is of class $C^1$, while for $F$ of class $C^k(k\geq 2)$ or real-analytic, $s\mapsto u(s)$ is of class $C^{k-1}$, respectively real-analytic.
\end{theo}

\begin{defi}
Let $\mathcal{A}$ be a set formed by the sequence pair $(\{\alpha_k\}, \{\beta_k\})$ satisfy the following conditions:\\
$(i)$ There exists some $k_0\in\mathbb{N}$, such that $\alpha_{k_0}\neq \beta_{k}, \forall k\in N$;\\
$(ii)$ For any $k,l\in\mathbb{N}$, there exist some $k_1, k_2, k_3, k_4, k_5, k_6$ such that
\[\alpha_k+\alpha_l=\alpha_{k_1}, \quad  \alpha_k-\alpha_l=\alpha_{k_2},\]
\[\beta_k+\beta_l=\alpha_{k_3}, \quad \beta_k-\beta_l=\alpha_{k_4},\]
\[\alpha_k+\beta_l=\beta_{k_5}, \quad \alpha_k-\beta_l=\beta_{k_6}.\]
Let $E^n$ be a subspace of $C^{n,\alpha}(\mathbb{R})\cap \mathcal{B}^2(\mathbb{R})/L_0$ which consists of those functions $f$ such that
\[f=\overline{f}+\sum_{k=0}a_k\cos(\alpha_k x)+\sum_{k=0}b_k\sin(\beta_k x),\]
where $(\{\alpha_k\}, \{\beta_k\})\in \mathcal{A}$. Moreover, we denote by $E^n_0$ a subspace of $E^n$ which consists of those functions $f$ with $\overline{f}=0$.
\end{defi}

\begin{rema}
Obviously, if $\alpha_k=Ck$ and $\beta_k=0$ with arbitrary constant $C>0$, then $(\{\alpha_k\}, \{\beta_k\})\in \mathcal{A}$. In this case, we see that $f$ is a periodic even function. By the above definition, we see that if $f\in C^{n,\alpha}(\mathbb{R})$ and $f$ is a periodic even function, then $f\in E^n$.
\end{rema}

\begin{rema}
The condition $(ii)$ of $(\{\alpha_k\}, \{\beta_k\})\in \mathcal{A}$ ensures that if $f\in E^n$ then $f^2\in E^n$ and $(f')^2\in E^{n-1}$.
\end{rema}

The following remarks give some examples about the sequence $(\{\alpha_k\}, \{\beta_k\})\in \mathcal{A}$ such that $\beta_k\neq 0$.

\begin{rema}
Let  $\alpha_k=2k$, $\beta_k=2k+1$, one can check that $(\{\alpha_k\}, \{\beta_k\})\in \mathcal{A}$, in this example $f$ is a periodic function.
\end{rema}

\begin{rema}\label{rema4.6}
Let $\alpha_{0}=\beta_0=0$. Proceed to define $\alpha_k$ and $\beta_k$ as follows: \\
Assuming that $k\geq 0$, $\alpha_k$ and $\beta_k$ are defined. Let
$$\alpha_{k+1}=\min_{n,l\in \mathbb{Z}}(2nA+2l B) \quad \textit{and} \quad \beta_{k+1}=\min_{n,l\in\mathbb{Z}}[(2n+1)A+(2l+1)B],$$
for some positive constant $A, B$ such that $\alpha_{k+1}\geq \alpha_k$ and $\beta_{k+1}\geq \beta_k$. If $A\in \mathbb{Q}^+, B\notin \mathbb{Q}^+$ (or $A\notin \mathbb{Q}^+, B\in \mathbb{Q}^+$), then one can check that $(\{\alpha_k\}, \{\beta_k\})\in \mathcal{A}$.
\end{rema}

Now we are going to utilize Theorem \ref{th1} to prove the existence of solutions to (\ref{4.2}). Let
\[X=\mathbb{R}\times E^n_0, \quad  Y=E^{n-1},\]
with $n\geq 1$ is an integer and $\alpha\in(0,1)$. The equation (\ref{4.2}) is equivalent to $F(\lambda,(\mu,w))=0$ where
\begin{align}\label{4.3}
F(\lambda,(\mu,w))&=\bigg\{\lambda+\gamma\bigg(\mathcal{G}_h(\frac{w^2}{2})-w-w\mathcal{G}_h(w)\bigg)\bigg\}^2-(\lambda^2+\mu-2gw)(w'^2+\mathcal{G}_h(w)^2+2\mathcal{G}_h(w)+1)\\
\nonumber&=\gamma^2\bigg(\mathcal{G}_h(\frac{w^2}{2})-w-w\mathcal{G}_h(w)\bigg)^2+2\lambda\gamma\bigg(\mathcal{G}_h(\frac{w^2}{2})-w-w\mathcal{G}_h(w)\bigg)\\
\nonumber&+(2gw-\mu)(w'^2+\mathcal{G}_h(w)^2+2\mathcal{G}_h(w)+1)-\lambda^2(w'^2+\mathcal{G}_h(w)^2+2\mathcal{G}_h(w)+1).
\end{align}
Obviously, $F\in C(\mathbb{R}\times X, Y)$ and $F(\lambda,(0,0))=0$ for all $\lambda\in \mathbb{R}$.
By virtue of (\ref{4.3}), we deduce that $\partial_\lambda F$, $\partial_{(\mu,w)}F$ exist and
\begin{align}
<\partial_{(\mu,w)}F(\lambda,(0,0)),(\nu,f)>=\lim_{t\rightarrow 0}\frac{F(\lambda,(t\nu,tf))-F(\lambda,(0,0))}{t}\\
\nonumber=2[(g-\lambda\gamma)f-\lambda^2\mathcal{G}_h(f)]-\nu,
\end{align}
where $(\nu,f)\in X$. Note that $f\in E_0$. Using the representation (\ref{2.1}), we obtain
\begin{multline}\label{4.5}
<\partial_{(\mu,w)}F(\lambda,(0,0)),(\nu,f)>=2(\sum^{\infty}_{k=1}a_k[(g-\lambda\gamma)-\lambda^2\alpha_k\coth(\alpha_kh))]\cos(\alpha_kx)\\
+2(\sum^{\infty}_{k=1}a_k[(g-\lambda\gamma)-\lambda^2\beta_k\coth(\beta_kh))]\sin(\beta_kx)-\nu.
\end{multline}
From the above equality, we see that if
\[\lambda^2\alpha_k\coth(\alpha_kh)\neq g-\lambda\gamma  \quad \text{and}\quad \lambda^2\beta_k\coth(\beta_kh)\neq g-\lambda\gamma ,\]
then $<\partial_{(\mu,w)}F(\lambda,(0,0)),(\nu,f)>=0$ if and only if $(\nu,f)=0$. Hence $\mathcal{N}(\partial_{(\mu,w)}F(\lambda,(0,0)))\neq{0}$ if and only if
\begin{align}
\lambda^2\alpha_k\coth(\alpha_kh)= g-\lambda\gamma  \quad \text{or}\quad \lambda^2\beta_k\coth(\beta_kh)= g-\lambda\gamma.
\end{align}
Since $f\in E_0$, it follows that there exists some $k_0$ such that $\alpha_{k_0}\neq \beta_{k_0}$. For this $k_0$, let us choose $\lambda^*$ such that
\[(\lambda^*)^2\alpha_{k_0}\coth(\alpha_{k_0}h)= g-\lambda^*\gamma.\]
By virtue of (\ref{4.5}), we deduce that the null space $\mathcal{N}(\partial_{(\mu,w)}F(\lambda^*,(0,0)))$ is one-dimensional and generated by $(0,\cos(\alpha_{k_0}x))\in X$. Moreover, $\mathcal{R}(\partial_{(\mu,w)}F(\lambda^*,(0,0)))$ is the closed subspace of $Y$ formed by the functions $f$ such that
\[\lim_{x\rightarrow 0}\frac{1}{x}\int^x_{-x}f(y)\cos(\alpha_{k_0} y)dy=0,\]
which implies that $Y/\mathcal{R}(\partial_{(\mu,w)}F(\lambda^*,(0,0)))$ is one-dimensional and generated by $\cos(\alpha_k x)$. Taking advantage of (\ref{4.5}), we have
\begin{align}
\partial^2_{\lambda,(\mu,w)}F(\lambda^*,(0,0))(1,(0,\cos(\alpha_{k_0} x)))&=2(-\gamma-2\lambda^*\alpha_{k_0}\coth(\alpha_{k_0} h))\cos(\alpha_{k_0} x)\\
\nonumber&=-2\lambda^*(\alpha_{k_0}\coth(\alpha_{k_0} h)+\frac{g}{(\lambda^*)^2})\cos(\alpha_{k_0} x).
\end{align}
Since $\cos(\alpha_{k_0} x)\notin \mathcal{R}(\partial_{(\mu,w)}F(\lambda^*,(0,0)))$ and $\lambda^*(\alpha_{k_0}\coth(\alpha_{k_0} h)+\frac{g}{(\lambda^*)^2})\neq 0$, it follows that
\[\partial^2_{\lambda,(\mu,w)}F(\lambda^*,(0,0))(1,(0,\cos(\alpha_{k_0} x)))\notin \mathcal{R}(\partial_{(\mu,w)}F(\lambda^*,(0,0))).\]

On the other hands, choose $\lambda^*$ such that
\[(\lambda^*)^2\beta_{k_0}\coth(\beta_{k_0}h)= g-\lambda^*\gamma.\]
In this case, we verify that the null space $\mathcal{N}(\partial_{(\mu,w)}F(\lambda^*,(0,0)))$ is one-dimensional and generated by $(0,\sin(\beta_{k_0}x))\in X$. Moreover, $\mathcal{R}(\partial_{(\mu,w)}F(\lambda^*,(0,0)))$ is the closed subspace of $Y$ formed by the functions $f$ such that
\[\lim_{x\rightarrow 0}\frac{1}{x}\int^x_{-x}f(y)\sin(\beta_{k_0} y)dy=0.\]
By the same token, one can deduce that $Y/\mathcal{R}(\partial_{(\mu,w)}F(\lambda^*,(0,0)))$ is one-dimensional and generated by $\sin(\beta_k x)$ and $\partial^2_{\lambda,(\mu,w)}F(\lambda^*,(0,0))(1,(0,\sin(\beta_{k_0}x))\notin \mathcal{R}(\partial_{(\mu,w)}F(\lambda^*,(0,0)))$.

Therefore, we obtain four bifurcation values
\[\lambda_1=-\frac{\gamma tanh(\alpha_{k_0}h)}{2\alpha_{k_0}}+\sqrt{\frac{\gamma^2 tanh^2(\alpha_{k_0}h)}{4\alpha^2_{k_0}}+g\frac{\tanh(\alpha_{k_0} h)}{\alpha_{k_0}}},\]
\[\lambda_2=-\frac{\gamma tanh(\alpha_{k_0}h)}{2\alpha_{k_0}}-\sqrt{\frac{\gamma^2 tanh^2(\alpha_{k_0}h)}{4\alpha^2_{k_0}}+g\frac{\tanh(\alpha_{k_0} h)}{\alpha_{k_0}}},\]
\[\lambda_3=-\frac{\gamma tanh(\beta_{k_0}h)}{2\beta_{k_0}}+\sqrt{\frac{\gamma^2 tanh^2(\beta_{k_0}h)}{4\beta^2_{k_0}}+g\frac{\tanh(\beta_{k_0} h)}{\beta_{k_0}}},\]
\[\lambda_4=-\frac{\gamma tanh(\beta_{k_0}h)}{2\beta_{k_0}}-\sqrt{\frac{\gamma^2 tanh^2(\beta_{k_0}h)}{4\beta^2_{k_0}}+g\frac{\tanh(\beta_{k_0} h)}{\beta_{k_0}}}.\]

Taking advantage of Theorem \ref{th1}, we prove the existence of solution with small amplitude.

Now we will give some example to explain the significance of our result.

\begin{rema}
Let $\alpha_k=k$ and $\beta_k=0$, then our result reduces to the result obtained in \cite{Constantin2011v}.
\end{rema}

\begin{rema}
Let $\alpha_k=2k$ and $\beta_k=2k+1$, we see that $(\{\alpha_k\},\{\beta_k\})\in\mathcal{A}$. In this case, if we choose $\lambda$ such that
\[\lambda=-\frac{\gamma tanh(h)}{2}-\sqrt{\frac{\gamma^2 tanh^2(h)}{4}+g\tanh(h)},\]
we have the local bifurcation curve
\[\{(\lambda(s),(o(s),s\sin(x)+o(s))):|s|<\ep\}\subset X.\]
This bifurcation value $\lambda$ is the same one as obtained in \cite{Constantin2011v} and the solution is also periodic. However, the solution obtained in \cite{Constantin2011v} is an even function. The solution constructed in our result is not an even function. So we obtain an non-uniqueness result corresponding to the same bifurcation value. Since
\[\lambda=\frac{m}{h}-\frac{\gamma h}{2},\]
it follows that for any given $m$, $h$, and the vorticity $\gamma$, there is at least two solutions that satisfy $(1.1)$. See Figure 1.
\end{rema}

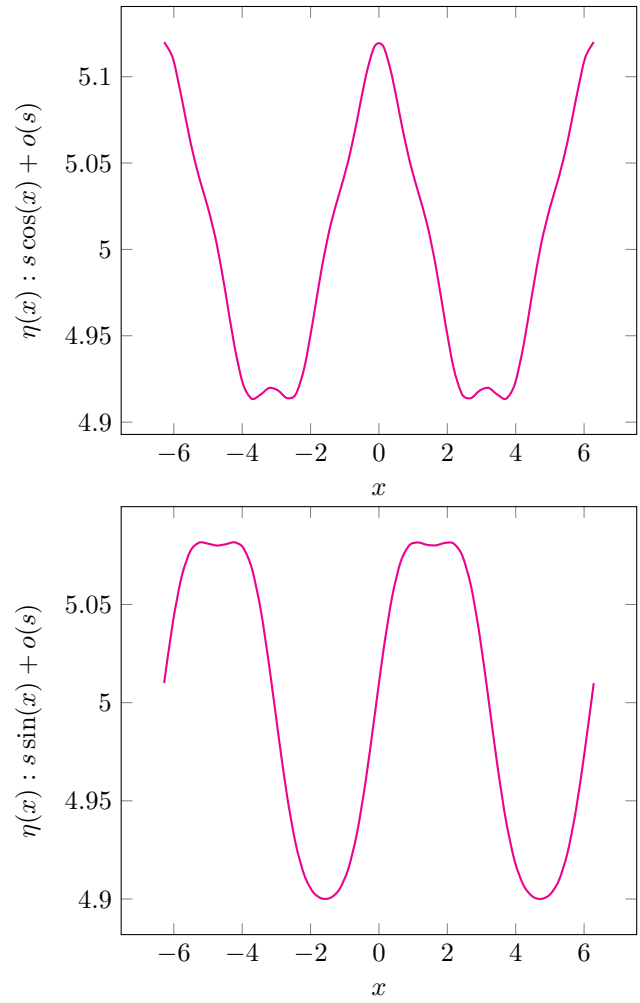
\begin{figure}
\centering
\begin{tikzpicture}
    \tikzset{elegant/.style={smooth,thick,samples=50,magenta}}
    \begin{axis}[ylabel=$\eta(x):s\cos(x)+o(s)$,
                 xlabel=$x$]
        \addplot[elegant,domain=-2*\num:2*\num]{1/10*cos(deg(x))+1/100*cos(2*deg(x))+1/100*cos(4*deg(x))+5};
    \end{axis}
\end{tikzpicture}
\begin{tikzpicture}
    \tikzset{elegant/.style={smooth,thick,samples=50,magenta}}
    \begin{axis}[ylabel=$\eta(x):s\sin(x)+o(s)$,
                 xlabel=$x$]
        \addplot[elegant,domain=-2*\num:2*\num]{1/10*sin(deg(x))+1/100*cos(2*deg(x))+1/100*sin(3*deg(x))+5};
    \end{axis}
\end{tikzpicture}
\caption{two different solutions at the same bifurcation point}\label{Fig1}
\end{figure}

\begin{rema}
Let $\alpha_k$ and $\beta_k$ defined as in Remark \ref{rema4.6} with $A=1$ and $B=\sqrt 5$, we can construct a solution that is almost periodic but not periodic (See Figure 2 below).
Thus, our result generalizes and covers the recent result obtained in \cite{Constantin2011v}.
\end{rema}

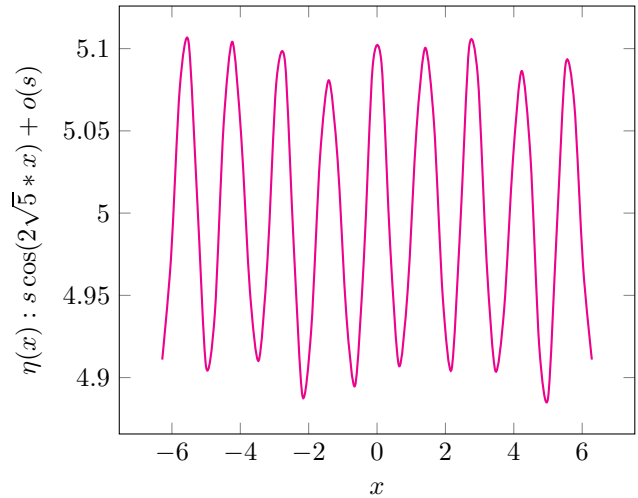
\begin{figure}
\centering
\begin{tikzpicture}
    \tikzset{elegant/.style={smooth,thick,samples=50,magenta}}
    \begin{axis}[ylabel=$\eta(x):s\cos(2\sqrt 5*x)+o(s)$,
                 xlabel=$x$]
        \addplot[elegant,domain=-2*\num:2*\num]{1/10*cos(2*sqrt(5)*deg(x))+1/100*sin(deg(x))+1/100*cos((4*sqrt(5)-2)*deg(x))+5};
    \end{axis}
\end{tikzpicture}
\caption{An almost periodic solution that is not a periodic function}\label{Fig2}
\end{figure}

\begin{rema}
The solution constructed in our result may contains the stagnation points. The stagnation points occur if only if
\[\lambda(\lambda+\gamma h)\leq 0.\]
Since the argument to show that the above inequality holds true is similar to  that of \cite{Constantin2011v}, we omit the details here.
\end{rema}

 {\bf Acknowledgements}. This work was
partially supported by the National Natural Science Foundation of China (No.11671407 and No.11701586), the Macao Science and Technology Development Fund (No. 098/2013/A3), and Guangdong Province of China Special Support Program (No. 8-2015),
and the key project of the Natural Science Foundation of Guangdong province (No. 2016A030311004).

\bibliographystyle{abbrv} 
\bibliography{myref}
\end{document}